\documentclass[10pt]{amsart}
\usepackage{amsfonts, amsthm, amssymb, amsmath, stmaryrd}
\usepackage{mathrsfs,array}
\usepackage{eucal,color,times,enumerate,accents}
\usepackage[all]{xy}
\usepackage{url}
\usepackage{enumitem}
\usepackage[pdftex,bookmarksnumbered,bookmarksopen]{hyperref}
\hypersetup{
  colorlinks,
  citecolor=red,
  linkcolor=blue,
  urlcolor=black}
  
    \newcommand\blfootnote[1]{%
  \begingroup
  \renewcommand\thefootnote{}\footnote{#1}%
  \addtocounter{footnote}{-1}%
  \endgroup
}

\newcommand{\Gad}{G^{\text{ad}}}
\newcommand{\Kad}{K^{\text{ad}}}
  
  \newcommand{\Addresses}{{
  \bigskip
  \footnotesize

\textsc{London School of Geometry and Number Theory, UCL, Department of Mathematics, Gower street, WC1E 6BT, London, UK}\par\nopagebreak
  \textit{E-mail address}, G.~Baldi: \texttt{gregorio.baldi.16@ucl.ac.uk}

}}

\usepackage{tikz-cd}
\usepackage{bbm}

\usetikzlibrary{positioning}
\usepackage{dsfont}
\usepackage{tikz}
\usetikzlibrary{matrix,arrows,decorations.pathmorphing,positioning}
\usepackage{calligra}
\usepackage[colorinlistoftodos]{todonotes}
\usepackage{xy}
\xyoption{all}

\DeclareMathOperator{\GSp2g}{GSp_{2g}}

\usepackage{tikz}
\usetikzlibrary{matrix,arrows,decorations.pathmorphing}

\usepackage{leftidx}

\input xy
\xyoption{all}

\usepackage[bottom=3.2cm, top=2.8cm, left=2.6cm, right=2.6cm]{geometry}

\usepackage{tikz-cd}
\theoremstyle{plain}
\newtheorem{thm}{Theorem}[section]

\newtheorem{prop}[thm]{Proposition}
\newtheorem{cor}[thm]{Corollary}
\theoremstyle{definition}
\newtheorem{defi}[thm]{Definition}

\newtheorem*{rmk}{Remark}

\theoremstyle{remark}

\numberwithin{equation}{subsection}

\newcommand{\pitop}{\pi_1^{\operatorname{top}}}
\newcommand{\piet}{\pi_1^{\operatorname{et}}}

\newcommand{\Gder}{G^{\text{der}}}

\DeclareMathOperator{\Sh}{Sh}

\DeclareMathOperator{\Hom}{Hom}

\DeclareMathOperator{\Res}{Res}

\DeclareMathOperator{\im}{Im}

\DeclareMathOperator{\Gl}{GL}

\newcommand{\Gm}{\mathbb{G}_m}

\newcommand{\MT}{\operatorname{MT}}

\newcommand{\Gsc}{G^{\text{sc}}}
\newcommand{\Xad}{X^{\text{ad}}}

\newcommand{\Z}{\mathbb{Z}}
\newcommand{\Q}{\mathbb{Q}}

\newcommand{\R}{\mathbb{R}}

\newcommand{\A}{\mathbb{A}}

\newcommand{\F}{\mathcal{F}}

\newcommand{\C}{\mathbb{C}}
\newcommand{\DT}{\mathbb{S}}

\newcommand{\Ff}{\mathbb{F}}

\makeatletter
\def\subtitle#1{\gdef\@subtitle{#1}}
\def\@subtitle{}
\makeatother

\begin{document}

\title{On the geometric Mumford-Tate conjecture for subvarieties of Shimura varieties}\blfootnote{\emph{Date}. May 22, 2019.}
\author{Gregorio Baldi}

\begin{abstract}
We study the image of $\ell$-adic representations attached to subvarieties of Shimura varieties $\Sh_K(G,X)$ that are not contained in a smaller Shimura subvariety and have \emph{no isotrivial components}. We show that, for $\ell$ large enough (depending on the Shimura datum $(G,X)$ and the subvariety), such image contains the $\Z_\ell$-points coming from the simply connected cover of the derived subgroup of $G$. This can be regarded as a geometric version of the integral $\ell$-adic Mumford-Tate conjecture.
\end{abstract}

\maketitle

\section{Introduction}
\subsection{Geometric and \texorpdfstring{$\ell$}{l}-adic monodromy}
Let $\DT$ denote the real torus $\Res_{\C / \R} (\Gm)$. A Shimura datum is a pair $(G,X)$ where $G$ is a reductive $\Q$-algebraic group and $X$ a $G(\R)$-orbit in the set of morphisms of $\R$-algebraic groups $\Hom(\DT, G_\R)$, satisfying the Shimura-Deligne axioms (\cite[Conditions 2.1.1(1-3)]{deligneshimura}). These axioms imply that the connected components of $X$ are hermitian symmetric domains and that faithful representations of $G$ induce variations of polarizable $\Q$-Hodge structures on $X$. Write $\A_f$ for the topological ring of finite adeles, and let $K$ be a neat compact open subgroup of $G(\A_f)$. Write
\begin{displaymath}
\Sh_K(G,X) := G(\Q)  \backslash \left ( X \times G(\mathbb{A}_f) / K \right).
\end{displaymath}
Let $X^+$ be a connected component of $X$ and $G(\Q)^+$ be the stabiliser of $X^+$ in $G(\Q)$. The above double coset set is a disjoint union of quotients of $X^+$ by the arithmetic groups $\Gamma_g:=G(\Q)^+ \cap gKg^{-1}$, where $g$ runs through a set of representatives for the finite double coset set $G(\Q)^+ \backslash G(\A_f )/K$. Baily and Borel \cite{MR0216035} proved that $\Sh_K(G,X)$ has a structure of a quasiprojective complex algebraic variety. 

For the length of this note we fix a faithful rational representation $G \subset \Gl( V)$ and an integral structure $V_\Z\subset V$ such that the image of $\Gamma := G(\Q)^+\cap K$ in $\Gl(V)$ is contained in $\Gl (V_\Z)$.

Let $S$ be the the image of $X^+ \times \{1\}$ in $ \Sh_K(G,X)$. Such varieties will be called \emph{connected components of Shimura varieties}. Write $\piet (S)$ for the étale fundamental group of $S$, with respect to some base point. One can attach to $S$ (and its subvarieties) an adelic representation
\begin{displaymath}
\piet (S)\to K \subset G(\A_f),
\end{displaymath}
corresponding to the tower of étale covering of $\Sh_K (G,X)$ indexed by open subgroups of $K$ (see for example \cite[Section 4]{MR3576114} and \cite[Section 2]{generalisedconj}).

Let $C \subset S$ be a smooth irreducible complex subvariety. We define the $\ell$-\emph{adic monodromy of C}, denoted by $\Pi_C^\ell$, as the image of 
\begin{displaymath}
\piet (C_\C) \to \piet (S_\C)\to K \xrightarrow{\pi_\ell} G(\Q_\ell),
\end{displaymath}
where $\pi_\ell : G(\A_f) \to G(\Q_\ell)$ denotes the projection to the $\ell$-th component. See Section \ref{comm} for a more detailed description of $\Pi_C^\ell$.

We want to prove that, for $\ell$-large enough, the $\ell$-adic monodromy of subvarieties $C \subset S$, satisfying two conditions we introduce in the next subsection, is \emph{large}.
\subsection{Main Theorem}
To study the geometric and $\ell$-adic monodromy of subvarieties of $S$ we may assume they are not contained in any smaller Shimura subvariety of $S$ (Theorem \ref{monodromyargument} will explain what happens if this is not the case). This is made precise in the next definition.
\begin{defi}
Let $S$ be a connected component of a Shimura variety $\Sh_K(G,X)$ and $C$ be an irreducible smooth complex subvariety of $S$. We say that $C$ is \emph{Hodge generic} if there exists a point $c\in C$  whose corresponding morphism of real algebraic groups $\DT \to G_\R$ does not factor through $H_\R \subset G_\R$ for any rational algebraic subgroup $H \varsubsetneq G$.
\end{defi}

Let $\Gad$ be the adjoint group of $G$ and $\Xad$ be the $\Gad(\R)$-orbit in $\Hom(\DT,\Gad)$ that contains the image of $X$ in $\Hom(\DT,\Gad)$. Choose a compact open subgroup $\Kad \subset \Gad(\A_f )$ containing the image of $K$. Then $(\Gad,\Xad)$ is a Shimura datum and we call $\Sh_{\Kad}(\Gad, \Xad)$ the \emph{adjoint Shimura variety} associated to $\Sh_K(G,X)$. The last ingredient needed to state our main theorem is the following.
\begin{defi}\label{nonisotrivial}
Let $S$ be a connected component of a Shimura variety $\Sh_K(G,X)$ and $C$ be an irreducible smooth complex subvariety of $S$. We say that $C$ has \emph{no isotrivial components} if, for every decomposition $(\Gad,\Xad, \Kad) \sim (G_1,X_1, K_1)\times \dots \times(G_n, X_n, K_n)$, the image of $C\to S \to \Sh_{K_i}(G_i, X_i)$ has dimension $>0$ for all $i=1,\dots, n$.
\end{defi}

Set $G(\Z_\ell)= G(\Q_\ell) \cap \Gl(V_\Z \otimes \Z_\ell)$. Let $\Gder$ be the derived subgroup of $G$ and $\lambda :\Gsc\to \Gder$ be the simply connected cover of $\Gder$. We write $G(\Z_\ell)^+$ for the subgroup of $G(\Z_\ell)$ given by the image of $\lambda_{\Z_\ell} : \Gsc(\Z_\ell)\to \Gder(\Z_\ell)$. It is an open subgroup of $\Gder(\Z_\ell)$ of index bounded independently of $\ell$. We now state the main theorem of this work.

\begin{thm}\label{main1}
Let $S$ be a connected component of a Shimura variety $\Sh_K(G,X)$ and $C$ be an irreducible smooth complex subvariety of $S$ which is Hodge generic and has no isotrivial components. For all $\ell$ big enough (depending only on $(G,X,K)$ and $C$), we have that $G(\Z_\ell)^+ \subset \Pi_C^\ell$. 
\end{thm}

\begin{rmk}In general one cannot expect the equality $ \Pi_C^\ell = G(\Z_\ell)$, even for large $\ell$. Indeed this may only happen when the Hodge structure associated to $C$, $h_C: \DT \to G_\R$, is \emph{Hodge maximal}, i.e. there is no non-trivial isogeny of connected $\Q$-groups $H\to G$ such that $h_C$ lifts to a homomorphism $h_C: \DT \to H_{\R} \to G_{\R}$ (see \cite[Definition 11.1]{sarremotivigalois} and \cite[Definition 2.1]{cadoretmoonen}). Of course if $G$ is simply connected, there are no such isogenies. Theorem \ref{main1}, when $G$ is semisimple simply connected, shows indeed that $G(\Z_\ell)^+ = G(\Z_\ell)=\Pi_C^\ell$, for all but finitely many primes $\ell$.
\end{rmk}

\subsection{Main Theorem for Shimura varieties of Hodge type}
Let $g$ be a positive natural number. The prototype for all Shimura varieties is the Siegel moduli space of principally polarized (complex) abelian varieties of dimension $g$ with a level structure, in which case $(G,X)=(\GSp2g, \mathcal{H}_g^\pm)$. For more details see \cite[Section 4]{delignetravaux}. We end the introduction reformulating Theorem \ref{main1} for subvarieties of Shimura varieties parametrising abelian varieties. We remark that the notion of having no isotrivial components for $C\subset \Sh_K (\GSp2g, \mathcal{H}_g^\pm)$ we proposed in the $(G,X)$-language (Definition \ref{nonisotrivial}) does not translate as parametrising a family of abelian varieties $A \to C$ with no non-zero constant factors (not even after passing to a finite cover of $C$). This is discussed in \cite[Scholie (page 13)]{andremt} and \cite[Section 4.5]{moonenlinearity}. Here we say that a $g$-dimensional abelian scheme $A\to C$ is \emph{with no isotrivial components} if it corresponds to a subvariety with no isotrivial components (as in Definition \ref{nonisotrivial}) of the Shimura variety associated to $(\GSp2g, \mathcal{H}_g^\pm)$. We have
\begin{cor}\label{cor}
Let $C$ be a complex smooth irreducible variety and $\eta$ its geometric generic point. Let $A \to C$ be a $g$-dimensional abelian scheme with no isotrivial components, write $M$ for the Mumford-Tate group of $A_{\eta}$ and $T_\ell (A_\eta)$ for the $\ell$-adic Tate module of $A_\eta$. For all $\ell$ large enough, the image of the map
\begin{displaymath}
\piet (C_\C) \to \Gl (T_\ell (A_\eta))
\end{displaymath}
contains $M(\Z_\ell)^+$.
\end{cor}
In proving the corollary, we may replace $C$ with a finite cover and the family $A \to C$ with an isogenous abelian variety. Hence we may assume that the family $A \to C$ gives rise to a subvariety of $\Sh_{K(3)}(\GSp2g, \mathcal{H}_g^\pm)$, where $K(3)\subset \GSp2g (\widehat{\Z})$ is the subgroup of elements that reduce to the identity modulo $3$. By assumption such subvariety is with no isotrivial components (in the sense of Definition \ref{nonisotrivial}). Since such family is contained in a Shimura subvariety whose defining group is the Mumford-Tate group of the abelian variety $A_\eta$, the corollary follows from Theorem \ref{main1} applied to the smooth locus of the image of $C$ in $\Sh_{K(3)}(\GSp2g, \mathcal{H}_g^\pm)$.

\begin{rmk}
We briefly explain why Corollary \ref{cor} can be thought as a geometric analogue of the \emph{integral Mumford-Tate} conjecture for abelian varieties (see \cite[Conjectures 10.3, 10.4 ,10.5]{sarremotivigalois} and \cite{cadoretmoonen}). Recall that, thanks to the work of Borovoi, Deligne, Milne and Milne-Shih, among others, Shimura varieties admit canonical models over number fields. Definition \ref{nonisotrivial} implies, among other things, that $C$ has positive dimension. Let us make the analogy with the case of an abelian variety $B$ over a number field $L$. Let $C$ be the spectrum of $L$. We still have a continuous morphism 
\begin{displaymath}
\rho_\ell : \piet (C) \to \Gl (T_\ell (B))
\end{displaymath} 
describing the action of the absolute Galois group of $L$ on the torsion points of $B$. The Mumford-Tate conjecture predicts that the image of $\rho_\ell$ is open in the $\Z_\ell$-points of the Mumford-Tate group $M$ of $B$ and its integral refinement that $\im (\rho_\ell)$ is as large as possible and its index in $M(\Z_\ell)$ can be bounded independently of $\ell$. 

Whenever the Shimura variety $\Sh_K(G,X)$ has an interpretation in terms of moduli space of motives, the same remark applies for Theorem \ref{main1}. We refer the reader to \cite[Section 4]{MR546619} for a discussion about when Shimura varieties are expected to have such interpretation and when they do not.
\end{rmk}

\subsection*{Conventions}
To simplify the notation, when writing the topological and étale fundamental groups, we omit the base-points.
\subsection*{Outline of paper}
In Section \ref{monodromy} we recall how to compute the geometric monodromy of subvarieties of (connected components of) Shimura varieties and explain the importance of Definition \ref{nonisotrivial}. In Section \ref{lastsection} we prove Theorem \ref{main1}, combining the results of Section \ref{monodromy} with a theorem of Nori.

\subsection*{Acknowledgements}
We thank Andrei Yafaev for his constructive suggestions. We also thank the anonymous referees whose comments simplified the proof of the main result and improved the exposition (spotting several inaccuracies).

\section{Monodromy of subvarieties, after Deligne, Andr\'{e} and Moonen}\label{monodromy}
In this section we recall how to produce, starting form a subvariety $C\subset \Sh_K(G,X)$, a Shimura subvariety of $\Sh_K(G,X)$ containing $C$ and such that $C$ becomes Hodge generic and has no isotrivial components in such Shimura subvariety. As a general reference for the theory of Shimura varieties we refer the reader to Deligne's works \cite{delignetravaux} and \cite{deligneshimura}.

\subsection{Preliminaries on Variations of Hodge structures}
Let $S$ be a connected, smooth complex algebraic variety, and $\mathcal{V}=(\mathcal{V},\F,\mathcal{Q})$ a polarized variation of $\Q$-Hodge structure on $S$. Let $\lambda: \widetilde{S}\to S$ be the universal cover of $S$ and fix a trivialisation $\lambda^* \mathcal{V}\cong \widetilde{S}\times V$. 
\begin{defi}
Let $s\in S$. The \emph{Mumford-Tate group at s}, denoted by $\MT_s \subset \Gl (\mathcal{V}_s)$, is the smallest $\Q$-algebraic group $M$ such that the map 
\begin{displaymath}
h_s : \DT \longrightarrow \Gl (V_{s,\R})
\end{displaymath}
describing the Hodge-structure on $V_{s}$, factors trough $M_\R$. Choosing a point $\tilde{s}\in \lambda^{-1}(s)\subset \widetilde{S}$ we obtain an injective homomorphism $\MT_s \subset \Gl(V)$. When a point $t \in S$ is such that $\MT_t$ is abelian (hence a torus), we say that $t$ is a \emph{special point}.
\end{defi}

It is well known, see \cite[Proposition 7.5]{delignek3}, that there exists a countable union $\Sigma \subsetneq S$ of proper analytic subspaces of $S$ such that
\begin{itemize}
\item for $s\in S - \Sigma $, $\MT_s \subset \Gl(V)$ does not depend on $s$, nor on the choice of $\tilde{s}$. We call this group \emph{the generic Mumford-Tate group of} $\mathcal{V}$ and we simply write it as $G$;
\item for all $s$ and $\tilde{s}$ as above, with $s\in \Sigma$, $\MT_s$ is a proper subgroup of $G$, the generic Mumford-Tate group of $\mathcal{V}$.
\end{itemize}

From now on assume that $\mathcal{V}$ admits a $\Z$-structure and choose $s\in S- \Sigma$ and $\tilde{s}$ as above. From the local system underlying $\mathcal{V}$ we obtain a representation $\rho: \pitop(S(\C)) \to \Gl(V)$, where $\pitop$ denotes the \emph{topological} fundamental group.
\begin{defi}
We denote by $M_s$ the connected component of the identity of the Zariski closure of the image of $\rho$ and we call $M_s$ the \emph{(connected) monodromy group}. Since we fixed a trivialisation for $\lambda^* \mathcal{V}$, we have that $M_s \subset \Gl(V)$ and $M_s$ does not depend on the choice of $s$ and $\tilde{s}$. 
\end{defi}
Recall the following important theorem. For the proof see \cite[Proposition 7.5]{delignek3}, \cite[Proposition 2]{andremt}; see also \cite[Section 1]{moonenlinearity}.
\begin{thm}[Deligne, Andr\'{e}, Moonen]\label{fixedpart}
Let $s\in S-\Sigma$. We have:
\begin{itemize}
\item[] \emph{Normality}. $M_s$ is a normal subgroup of the derived group $G^{\text{der}}$;
\item[] \emph{Maximality}. Suppose $S$ contains a special point. Then $M_s=G^{\text{der}}$.
\end{itemize}
\end{thm}

An immediate application of Theorem \ref{fixedpart} is the following. Let $S$ be a connected component of a Shimura variety $\Sh_K(G,X)$. For simplicity assume that $G$ is a semi-simple algebraic group of adjoint type and that $G$ is the generic Mumford-Tate group on $X$.  Fixing a rational representation of $G$ in $\Gl (V)$, as we did in the introduction, we obtain a polarized variation of $\Q$-Hodge structure, denoted by $\mathcal{V}$, on the constant sheaf $V_{X^+}$. Moreover, since $\Gamma=G(\Q)^+\cap K$ acts freely on $X^+$, $\mathcal{V}$ descends to a variation of Hodge structures on $S$. To obtain a $\Z$-structure, and apply the previous theorem, we may choose $V_{\widehat{\Z}}$, a $K$-invariant lattice in $V\otimes \A_f$, and define $V_\Z$ as $V \cap V_{\widehat{\Z}}$. Let $c\in C$ be a Hodge-generic point on $C$, since $\pi_1(C,c)$ acts on $\mathcal{V}_{\Z,c}$, it acts on $V_\Z$. Let $\Pi$ be its image in $\Gl(V_\Z)$, it is a finitely generated group. We have

\begin{cor}\label{cormonodromy}
Let $C$ be a smooth irreducible Hodge generic subvariety of $S$ containing a special point $t$. We have $\Pi \subset \Gamma$ and both of them are Zariski dense in $G$.
\end{cor}
The next section explains what happens when $C$ does not contain a special point and it is not Hodge generic. More details can be found in \cite[Sections 2.9, 3.6 and 3.7]{moonenlinearity}.
\subsection{Monodromy of subvarieties with no isotrivial components}\label{moonen}
Let $C$ be an irreducible smooth complex subvariety of a Shimura variety $\Sh_K(G,X)$. Since the intersection of two Shimura subvarieties of $\Sh_K(G,X)$ is again a Shimura subvariety, there exists a unique smallest sub-Shimura variety $S_C \subset\Sh_K(G,X)$ containing $C$. By definition there exists a $\Q$-group $M$ such that $S_C$ is an irreducible component of the image of $X^+_M \times \eta K$ in $\Sh_K(G,X)$, for some $\eta \in G(\A_f)$, where $X_M^+$ is the restriction of $X^+$ to $M$. Moreover we may take $M$ to be the generic Mumford-Tate group on $C$ (recall that we fixed from the beginning a faithful representation of $G$). By construction, $C$ is Hodge generic in $S_C$.

Let $H$ the connected monodromy group associated to the polarized variation of $\Z$-Hodge structures $\mathcal{V}$ restricted to $C$. Theorem \ref{fixedpart} implies that $H$ is a normal subgroup of the derived subgroup of $M$ and, since $M$ is reductive, we can find a normal algebraic subgroup $H_2$ in $M$ such that $M$ is the almost direct product of $H$ and $H_2$. This induces a decomposition of the adjoint Shimura datum:
\begin{displaymath}
(M^{\operatorname{ad}}, X_{M^{\operatorname{ad}}})=(H^{\operatorname{ad}}, X_{H^{\operatorname{ad}}})\times (H^{\operatorname{ad}}_2, X_{H^{\operatorname{ad}}_2}).
\end{displaymath}
As in Corollary \ref{cormonodromy}, when $C$ contains a special point, $H$ is exactly the derived subgroup of $M$, $X_{M}$ is isomorphic to $X_{H}$ and $X_{H_2}$ is nothing but a point. But if $C$ does not contain a special point, Moonen, in \cite[Proposition 3.7]{moonenlinearity}, proves the following. 

\begin{prop}[Moonen]\label{propmoonen}
Let $\mathcal{C}$ be an irreducible component in the preimage of $C$ in $X_M$. The image of $\mathcal{C}$ under the projection $X_M \to X_{H_2}$ is a single point, say $y_2\in X_{H_2}$. We have that $C$ is contained in the image of $(Y_1\times \{y_2\}) \times \eta' K$ in $\Sh_K(G,X)$ for some connected component $Y_1 \subset X_{H} $ and a class $\eta' K \in G(\A_f)/K$.
\end{prop}

To summarise our discussion, from Corollary \ref{cormonodromy} and Proposition \ref{propmoonen} applied to $C$, we have
\begin{thm}\label{monodromyargument}
Let $C$ be an irreducible smooth complex Hodge generic subvariety of a Shimura variety $\Sh_{K_M}(M,X_M)$. There exists a unique sub-Shimura datum $(H, X_H)\hookrightarrow (M,X_M)$ such that
\begin{displaymath}
(M^{\operatorname{ad}}, X_{M^{\operatorname{ad}}})=(H^{\operatorname{ad}}, X_{H^{\operatorname{ad}}})\times (H^{\operatorname{ad}}_2, X_{H_2^{\operatorname{ad}}}),
\end{displaymath}
and \begin{itemize}
\item the projection of $C$ to $\Sh_{K_{H^{\operatorname{ad}}}} (H^{\operatorname{ad}}, X_{H^{\operatorname{ad}}})$, denoted by $\widetilde{C}$, is Hodge generic and has no isotrivial components;
\item the projection of $C$ to $\Sh_{K_{H_2^{\operatorname{ad}}}}(H^{\operatorname{ad}}_2, X_{H_2})$ is a single point.
\end{itemize}
Moreover
\begin{itemize}
\item the fundamental groups of $\widetilde{C}$ and of $\Sh_{K_{H^{\operatorname{ad}}}} (H^{\operatorname{ad}}, X_{H^{\operatorname{ad}}})$ are both Zariski dense in $H^{\operatorname{ad}}$.
\end{itemize}
\end{thm}
\section{\texorpdfstring{$\ell$}{l}-adic monodromy}\label{lastsection}
After recalling the interplay between the geometric and the $\ell$-adic monodromy of subvarieties of $S$, we eventually prove Theorem \ref{main1}.
\subsection{A commutative diagram}\label{comm}
Let $\mathcal{V}$ be a variation of polarized $\Z$-Hodge structures on $S$ (as explained in Section \ref{monodromy}) and let
\begin{displaymath}
\pitop (C(\C)) \to \Gamma \subset \Gl (V_\Z)
\end{displaymath}
be the monodromy representation of the induced variation on $C$. Let $\overline{\Gamma}$ be the closure of $\Gamma$ in $\Gl (V_\Z \otimes \widehat{\Z})$. The map $\pitop(C(\C)) \to\Gamma \to\overline{\Gamma}$ canonically factorises trough the profinite completion of $\pitop(C(\C))$, which is canonically isomorphic to $\piet(C_\C)$, the  \'{e}tale fundamental group of $C_\C$. Therefore we have a commutative diagram
\begin{center}
\begin{tikzpicture}[scale=1.5]
\node (A) at (0,1) {$\pitop(C(\C))$};
\node (B) at (2,1) {$\Gamma \subset \Gl (V_\Z)$};
\node (C) at (0,0) {$\piet(C_\C)$};
\node (D) at (2,0) {$\overline{\Gamma} \subset \Gl (V_\Z \otimes \widehat{\Z}) $};
\path[->,font=\scriptsize,>=angle 90]
(A) edge node[above]{} (B)
(A) edge node[right]{} (C)
(B) edge node[right]{} (D)
(C) edge node[above]{} (D);
\end{tikzpicture}.
\end{center}
Let $\pi_\ell: \Gl (V_\Z \otimes \widehat{\Z})\to \Gl (V_\Z \otimes \Z _\ell)$ be the projection to the $\ell$-th component. The image of the map
\begin{displaymath}
\piet(C_\C)\to \overline{\Gamma}\xrightarrow{ \pi_\ell} \Gl (V_\Z \otimes \Z_\ell)
\end{displaymath}
is precisely the $\ell$-adic monodromy of $C$, that we denoted in the introduction by $\Pi_C^\ell$.

\subsection{Proof of Theorem \ref{main1}}
First of all notice that, to prove Theorem \ref{main1}, the difference between $(\Gad,\Xad)$ and $(G,X)$ is irrelevant. Indeed the arithmetic data play no role in the problem. Therefore we may, and do, assume that $G$ is of adjoint type. Moreover, since the Theorem is about large enough primes $\ell$, we may ignore finitely many primes, and assume that $G$ is the generic fiber of a semisimple adjoint group scheme over $\Z$.

Before starting the proof we fix some notations. By $(-)_\ell$ we denote the reduction modulo-$\ell$ of subgroups of $G(\Z)$. Write $G(\Ff_\ell)$ (resp. $G(\Ff_\ell)^+$) for the reduction modulo-$\ell$ of $G(\Z)$ (resp. $G(\Z_\ell)^+$).

Let $C\subset S=\Gamma \backslash X^+$ be a subvariety as in the statement of Theorem \ref{main1} and let $\Pi\subset \Gamma$ be the image of the geometric monodromy representation
\begin{displaymath}
\pitop(C(\C)) \to \pitop(S(\C))\to\Gamma \subset \Gl(V_\Z).
\end{displaymath}
Since $C$ is Hodge generic and has no isotrivial components, Corollary \ref{cormonodromy} shows that $\Pi$ is Zariski dense in $G$.

We recall a theorem of Nori (\cite[Theorem 5.1]{nori}) about Zariski dense subgroups of semisimple groups.
\begin{thm}[Nori]\label{nori}
Let $H \subset {\Gl_n} / \Q$ be a semisimple group and $\Pi \leq H(\Q)$ be a discrete finitely generated Zariski-dense subgroup. Then for all sufficiently large prime numbers $\ell$ (depending only on $H$ and $\Pi$), the reduction modulo-$\ell$ of $\Pi$ contains $H(\Ff_\ell)^+$.
\end{thm}
Reducing modulo-$\ell$ is well defined for $\ell$ large enough: indeed, since $\Pi$ is finitely generated, there are only finitely many primes $\ell_1, \dots, \ell_k$ such that $\Pi$ belongs to
\begin{displaymath}
H(\Z[\ell_i^{-1}]_i):= H \cap \Gl_n(\Z[\ell_1^{-1},\dots, \ell_k^{-1}]),
\end{displaymath}
and the reduction mod $\ell$ is well-defined for all other primes.

Since the group $G$ is assumed to be semisimple and adjoint we may apply Theorem \ref{nori} to get, for $\ell$ large enough, a chain of inclusions
\begin{displaymath}
G(\Ff_\ell)^+ \subset \Pi_\ell \subset \Gamma_\ell \subset G(\Ff_\ell).
\end{displaymath}

We are left to lift such chain of inclusions to the $\Z_\ell$-points of $G$. Denote by $\alpha_\ell$ the reduction modulo-$\ell$ map:
\begin{displaymath}
\alpha_\ell : G(\Z_\ell)\to G(\Ff_\ell).
\end{displaymath}
It is a well known fact that, for $\ell$ large enough, if $G$ is a connected semisimple group, $\alpha_\ell$ is a \emph{Frattini cover}, i.e. $\alpha_\ell$ is surjective and $G(\Z_\ell)$ contains no strict subgroups mapping surjectively onto $G(\Ff_\ell)$ (equivalently, the kernel of $\alpha_\ell$ is contained in the Frattini subgroup of $G(\Z_\ell)$). A proof of this fact can be found in \cite[Lemma 16.4.5 (page 403)]{lubotzky}, see also \cite[Section 2.3]{MR3455865} and \cite[Proposition 7.3]{MR735226}. 

Nori, in \cite[Section 3]{nori} (in particular \cite[Remark 3.6]{nori}), shows also that, for $\ell$ large enough, $G(\Ff_\ell)^+$ can be identified with the subgroup of $G(\Ff_\ell )$ generated by its $\ell$-Sylow subgroups. Hence we have that the index $[G(\Ff_\ell):G(\Ff_\ell)^+]$, for $\ell$-large enough, is prime to $\ell$. From this, and the fact that $\alpha_\ell $ is Frattini, we deduce that the maps
\begin{displaymath}
 \alpha_\ell^{-1}(G(\Ff_\ell)^+)\to G(\Ff_\ell)^+ , \text{   and   } \alpha_\ell^{-1}(\Pi_\ell) \to \Pi_\ell
\end{displaymath}
are also Frattini covers (see also \cite[Section 2.3]{MR3455865}). Since $\Pi^\ell_C \subset \alpha_\ell^{-1}(\Pi_\ell) $ surjects onto $\Pi_\ell$ and $G(\Z_\ell)^+$ onto $G(\Ff_\ell)^+$, the inclusions
\begin{displaymath}
G(\Z_\ell)^+\subset \alpha_\ell^{-1}(G(\Ff_\ell)^+),  \text{   and   } \Pi^\ell_C \subset \alpha_\ell^{-1}(\Pi_\ell)
\end{displaymath}
are actually equalities. Eventually we have
\begin{displaymath}
G(\Z_\ell)^+\subset  \Pi^\ell_C,
\end{displaymath}
as desired. This ends the proof of Theorem \ref{main1}.

\bibliographystyle{alpha}
\bibliography{biblio.bib}

\Addresses

\end{document}